\documentclass[12pt]{article}
\usepackage{epsfig}
\usepackage{latexsym}
\usepackage{amsfonts}
\usepackage{amsmath}
\usepackage[psamsfonts]{amssymb}

\def\K{{\mathbf{K}}}

\def\R{{\mathbf{R}}}
\def\S{{\mathbf{S}}}

\def\P{{\mathbf{P}}}
\def\C{{\mathbf{C}}}

\newtheorem{thm}{Theorem}[section]
\newtheorem{example}[thm]{Example}
\newtheorem{prop}[thm]{Proposition}

\newtheorem{df}[thm]{Definition} 

\numberwithin{equation}{section}

\author{Alexandre Eremenko\thanks{Supported by NSF grant DMS-1351836.},
Andrei Gabrielov\thanks{Supported by NSF grant DMS-1161629.} and Vitaly Tarasov}
\title{Metrics with conic singularities and spherical polygons}
\begin{document}
\maketitle
\begin{abstract} A spherical $n$-gon is a bordered surface
homeomorphic to a closed disk, with $n$ distinguished
boundary points called corners,
equipped with a Riemannian
metric of constant curvature $1$, except at the corners, and such that
the boundary arcs between the corners are geodesic.
We discuss the problem of classification of these polygons
and enumerate them in the case that
two angles at the corners are not multiples of $\pi$.
The problem is equivalent to classification of some second order
linear differential
equations with regular singularities, with real
parameters and
unitary monodromy.


\end{abstract}

\noindent
\section{Introduction}\label{sec:Intro}

Let $S$ be a compact Riemann surface, and $\{ a_0,\ldots,a_{n-1}\}$
a finite set
of points on $S$. Let us consider a conformal Riemannian metric on $S$
of constant curvature
$K\in\{0,1,-1\}$ with conic singularities at the points $a_j$.
This means that in a local conformal coordinate $z$ the length element of the
metric is given by the formula $ds=\rho(z)|dz|$, where
$\rho$ is a solution of the differential equation
\begin{equation}\label{0}
\Delta\log\rho+K\rho^2=0\quad{\mbox{in}}
\quad S\backslash \{ a_0,\ldots,a_{n-1} \},
\end{equation}
and $\rho(z)\sim |z|^{\alpha_j-1}$, for the local coordinate $z$ which
is equal to $0$ at $a_j$. Here $\alpha_j>0,$ and $2\pi\alpha_j$ is
the total angle around the singularity $a_j$.

Alternatively, for every point in $S$, there exists a local coordinate
$z$ for which
$$ds=\frac{2\alpha|z|^{\alpha-1}|dz|}{1+K|z|^2},$$
where $\alpha>0$.
At the singular points $a_j$ we have $\alpha=\alpha_j$ while at
all other points, $\alpha=1$.

In 1890, the G\"ottingen Mathematical Society proposed the study of
equation (\ref{0}) as a competition topic, probably by
suggestion of H. A. Schwarz \cite{SG}.
E. Picard wrote several papers on the
subject, \cite{P1,P2,P}, see also \cite[Chap. 4]{P3}.
When $K<0$, the topic is closely related to the uniformization of orbifolds
\cite{SG,For,Po}. In the case of uniformization of orbifolds,
one is interested in the
angles $2\pi\alpha_j=2\pi/m_j$, where $m_j$ are positive integers.
The case $K\leq 0$ is quite well understood, but very little is known on
the case $K>0$.

McOwen \cite{MO} and
Troyanov \cite{T} studied the general question of existence and uniqueness
of such metrics with prescribed $a_j$, $\alpha_j$ and $K$.
Troyanov also considered
the case of non-constant curvature $K$.
One necessary condition that one has to impose on these data follows
from the Gauss--Bonnet theorem:
the quantity
\begin{equation}\label{curv}
\chi(S)+\sum_{j=0}^{n-1}(\alpha_j-1)\quad\mbox{has the same sign as}\quad K.
\end{equation}
Here $\chi$ is the Euler characteristic.
Indeed, this quantity multiplied by
$2\pi$ is the integral curvature of the smooth part of the surface.

It follows from the results of Picard, McOwen
and Troyanov, that for $K\leq 0$,
condition (\ref{curv})
is also sufficient for the existence of the metric with conic singularities
at arbitrary points $a_j$ and angles $2\pi\alpha_j$. The metric with
given $a_j$ and $\alpha_j$ is unique
when $K<0$, and unique up to a constant multiple when $K=0$.

In the case of positive curvature, the results are much less complete.
The result of Troyanov that applies to $K=1$ is the following:
\vspace{.1in}

{\em Let $S$ be a compact Riemann surface, $a_0,\ldots,a_{n-1}$
points
on $S$, and \newline
$\alpha_0,\ldots,\alpha_{n-1}$ po\-si\-ti\-ve num\-bers
sa\-tis\-fy\-ing
\begin{equation}\label{troy}
0<\chi(S)+\sum_{j=0}^{n-1}(\alpha_j-1)
<2\min\{1,\alpha_0,\ldots,\alpha_{n-1}\}.
\end{equation}
Then there exists a conformal metric of positive curvature $1$ on $S$
with conic singularities at $a_j$ and angles $2\pi\alpha_j$.}
\vspace{.1in}

F. Luo and G. Tian \cite{LT} proved that if the condition
$0<\alpha_j<1$ is satisfied for $0\leq j\leq n-1,$
then (\ref{troy}) is necessary and sufficient,
and the metric with given $a_j$ and $\alpha_j$ is unique.

In general, the second
inequality in (\ref{troy})
is not necessary.
In this paper, we only
consider the simplest case when $S$ is the sphere, so $\chi(S)=2$.
The problem of description and classification of conformal metrics
of curvature $1$ with conic singularities on the sphere has
applications to the study of certain surfaces of constant mean curvature
\cite{FKKRUY,D2,D}, and to several other questions of
geometry and physics \cite{HR,T}.

The so-called ``real case'' is interesting and important.
Suppose that all singularities belong to a circle on the sphere $S$,
and we only
consider the metrics which are symmetric with respect to this circle.
Then the circle splits $S$ into two symmetric halves, and each of them
is a {\em spherical polygon},
for which we state a formal definition:

\begin{df}\label{n-gon}{\rm {\em A spherical $n$-gon} is a closed disk
with $n$ distinguished boundary points $a_j$ called the corners,
equipped with
a conformal Riemannian metric of constant curvature
$1$ everywhere except the corners, and such that the sides (boundary arcs
between the corners) are geodesic. The metric has conic singularities at
the corners.}
\end{df}

\begin{example}\label{flat-n-gons}
Let us consider {\em flat} $n$-gons, which are defined
similarly, with $K=0$.
The necessary and sufficient
condition for the existence of a flat $n$-gon with angles $\pi\alpha_j$
is given by the Gauss--Bonnet theorem which in this case says that
$$\sum\alpha_j=n-2,$$
and the polygon with given angles and prescribed corners\footnote{By this we
mean that the images of the corners under a conformal map onto
a disk are prescribed.}
is unique up to a scaling factor. The simplest proof of these facts
is the Schwarz--Christoffel formula.
Thus our subject can be considered as a generalization of the
Schwarz--Christoffel formula to the case of positive curvature.
\end{example}
\vspace{.1in}

In \cite{UY,E}, all possibilities
for spherical triangles are completely described,
see also \cite{FKKRUY} where a minor error in \cite[Theorem 2]{E}
is corrected. In the case of triangles,
the metric is uniquely determined by the angles in the case
that none of the $\alpha_j$ is an integer.

The case when all $\alpha_j$ are integers, and $n$ is arbitrary,
is also well understood.
In this case, the line element of the metric has the global representation
$$ds=\frac{2|f'||dz|}{1+|f|^2},$$
where $f$ is a rational function. The singular points $a_j$ are the critical
points of $f$, $\alpha_j-1$ is the multiplicity
of the critical point at $a_j$, and $\alpha_j$ is the local degree
of $f$ at $a_j$.

Thus the problem with all integer $\alpha_j$ is equivalent to describing
rational functions with prescribed critical points \cite{G,S,EG0,EG,EG2,EGSV}.
Also, this problem is equivalent to diagonalization of commuting Hamiltonians
for the quantum Gaudin model by the Bethe ansatz method
\cite{Gaud,MTV,MTV3}.

Almost nothing is known in the case that
some of the $\alpha_j$ are not integers, the number of singularities
is greater than $3$, and the second
inequality in (\ref{troy}) is violated.

In this paper we investigate the case when two of the $\alpha_j$
are non-integer while the rest are integers,
with the emphasis on the real case. See also \cite{CWX} where some examples
of such metrics are given.

Our method can be considered a generalization of
the method in \cite{EG0,EG,EGSV} and Klein \cite{Klein},
who
classified circular triangles with arbitrary angles,
not necessarily geodesic.
A modern paper which uses Klein's approach to triangles is \cite{Y}.

The contents of the paper is the following.
In section \ref{DE-connection}
we recall the connection of the problem with linear differential
equations.
In section \ref{all-integers} we describe
what is known about the case when all $\alpha_j$ are integers,
with a special emphasis on the ``real case'' when all $a_j$ belong
to the real line.

In sections \ref{cond-angles}-\ref{deformation} we
study the case when two
of the $\alpha_j$ are not integers and the others are integers, and give
a complete classification of spherical polygons in this case.

\section{Connection with linear differential equations}\label{DE-connection}

Let $(S,ds)$ be the Riemann sphere equipped
with a metric with conic singularities.
Every smooth point of $S$ has a neighborhood which is isometric to a region on
the standard unit sphere $\S$; let $f$ be such an isometry.
Then $f$ has an analytic continuation along every path
in $S\backslash\{ a_0,\ldots,a_{n-1}\}$, and we obtain a multi-valued function
which is called the {\em developing map}. The monodromy of
$f$ consists of orientation-preserving isometries (rotations) of
$\S$, so the Schwarzian derivative
\begin{equation}
\label{schwarz}
F(z):=\frac{f^{\prime\prime\prime}}{f^\prime}-\frac{3}{2}\left(
\frac{f^{\prime\prime}}{f^\prime}\right)^2
\end{equation}
is a single valued function.

Developing map is completely characterized by the properties that
it has an analytic continuation along any curve in
$S\backslash\{ a_0,\ldots,a_{n-1}\}$,
has asymptotics $\sim c(z-a_j)^{\alpha_j}$ as $z\to a_j,$ where $c\neq 0,$
and has
$PSU(2)=SO(3)$ monodromy. It is possible that two such maps with the same
$a_j$ and $\alpha_j$ are related by
post-composition with a fractional-linear transformation.
The metrics arising from such maps will be called {\em equivalent}.
Following \cite{FKKRUY}, we say that the metric is reducible
if its monodromy group is commutative (which is equivalent to
all monodromy transformations having a common fixed point).
In the case of irreducible metrics, each equivalence class
contains only one metric. For reducible metrics, the equivalence
class is a
$1$-parametric family when the monodromy is non-trivial and
a $2$-parametric family when monodromy is trivial.

The asymptotic
behavior of $f$ at the singular
points $a_j$ implies that the only singularities of $F$ on
the sphere are double poles, so $F$ is a rational function,
and we obtain the Schwarz differential equation (\ref{schwarz}) for $f$.

It is well-known that the general solution
of the Schwarz differential equation
is a ratio of two linearly independent solutions of
the linear differential
equation
\begin{equation}
\label{linear}
y^{\prime\prime}+Py'+Qy=0,\quad f=y_1/y_0,
\end{equation}
where $$F=-P'-P^2/2+2Q.$$
For example one can take $P=0$, then $Q=F/2$.
Another convenient choice is to make all poles but one of $P$ and $Q$
simple.
When $n=3$, equation (\ref{linear})
is equivalent to the hypergeometric equation,
and when $n=4$ to the Heun equation \cite{R}.

The singular points $a_j$ of the metric
are the singular points of the equation
(\ref{linear}). These singular points are regular,
and to each point correspond
two exponents $\alpha_j^\prime>\alpha_j^{\prime\prime}$,
so that $\alpha_j=\alpha_j^\prime-\alpha_j^{\prime\prime}$.
If $\alpha_j$ is an integer for some $j$,
we have an additional condition
of the absence of logarithms in the formal solution
of (\ref{linear}) near $a_j$.

It is easy to write down the general form of a Fuchsian equation with
prescribed
singularities and prescribed exponents at the singularities.
By a fractional-linear change of the independent variable one
can achieve that one singular point, say $a_{n-1}$ is $\infty$.
Then by
chan\-ges of the va\-riable $y(z)\mapsto y(z)(z-a_j)^{\beta_j}$,
one can achieve that the smaller expo\-nent at each finite singular point
is $0$, see \cite{R}.
After a normalization, $n-3$ parameters remain,
the so-called accessory parameters \cite[Ch. IV, 7]{Golub}
\begin{equation}\label{ode}
w^{\prime\prime}+\sum_{j=0}^{n-2}\frac{1-\alpha_j}{z-a_j}w^\prime+
\frac{\alpha'\alpha^{\prime\prime}z^{n-3}+\lambda_{n-4}z^{n-4}+\ldots+
\lambda_0}{\prod_{j=0}^{n-2}(z-a_j)}w=0,
\end{equation}
where
\begin{equation}\label{cond}
\alpha_{n-1}=\alpha'-\alpha^{\prime\prime},\quad
\sum_{j=0}^{n-2}\alpha_j+\alpha^\prime+\alpha^{\prime\prime}=n-2.
\end{equation}
Here the exponents at the singular points are described
by the Riemann symbol
$$P\left\{\begin{array}{cccccc}a_0&a_1&a_2&\ldots&\infty&\\
0&0&0&\ldots&\alpha^{\prime\prime}&;z\\
\alpha_0&\alpha_1&\alpha_2&\ldots&\alpha^{\prime}&
\end{array}\right\}.$$
The first line lists the singularities, the second the smaller
exponents and the third the larger exponents.
So the angle at infinity is $\alpha_{n-1}=\alpha^\prime-\alpha^{\prime\prime}$.
The accessory parameters are $\lambda_0,\ldots,\lambda_{n-4}$.

To obtain a conformal metric
of curvature $1$, one has to choose these accessory parameters
in such a way that the monodromy group
of the equation is conjugate to a subgroup
of $PSU(2)$.

Solving (\ref{cond}),
we obtain
\begin{equation}
\label{alphaprime}
\alpha'=\frac{1}{2}\left(n-2+\alpha_{n-1}-\sum_{j=0}^{n-2}\alpha_j\right)
\end{equation}
and
$$\alpha^{\prime\prime}=\frac{1}{2}\left(n-2-\alpha_{n-1}-\sum_{j=0}^{n-2}
\alpha_j\right).
$$
The question of existence of a spherical metric
with given singularities
\newline
$a_0,a_1,\ldots,a_{n-2},\infty$
and given angles $2\pi\alpha_j,\; 0\leq j\leq n-1$
is equivalent to the following:
{\em When one can choose the accessory parameters $\lambda_j$
so that the projective monodromy
group of the equation (\ref{ode})
is conjugate to a subgroup of $PSU(2)$ ?}

The necessary condition (\ref{curv}) can be restated for the equation
(\ref{ode})
as
\begin{equation}
\label{alphapp}
\alpha^{\prime\prime}<0.
\end{equation}
For the case of four singularities, questions similar to our problem
were investigated in \cite{Klein2,Hilb1,Sm,Sm2}:
{\em when can one choose the accessory parameter so that the monodromy group
of the Heun equation preserves a circle?}
All these authors consider the problem under the assumption
\begin{equation}
\label{assumption}
0\leq\alpha_j<1,\quad\mbox{for}\quad 0\leq j\leq n-1.
\end{equation}
In the paper \cite{D} the problem of choosing the accessory
parameter so that the monodromy group is conjugate to a subgroup in $PSU(2)$ is discussed.
However all results of this
paper are also proved only under the assumption (\ref{assumption}).
Assumption (\ref{assumption}) seems to be essential for the methods
of Klein \cite{Klein2}, Hilb \cite{Hilb1}, Smirnov \cite{Sm,Sm2},
Dorfmeister and Schuster \cite{D,D2}.

\section{Case when all $\alpha_j$ are integers}\label{all-integers}

If all $\alpha_j$ are integers, the developing map $f$ is a rational
function, and the metric of curvature $1$
with conic singularities can be globally described
as the pull-back of the spherical metric via $f$,
that is
$$ds=\frac{2|f'||dz|}{1+|f|^2}.$$
The singular points of the metric are critical points of $f$, and
$\alpha_j-1$ are multiplicities of these critical points.

The following results are known for this case.

First of all, the sum of $\alpha_j-1$
must be even: if $d$ is the degree of $f$, then
\begin{equation}\label{1}
2+\sum_{j=1}^n(\alpha_j-1)=2d.
\end{equation}
This is stronger than the
general necessary condition (\ref{curv}).

Second,
\begin{equation}\label{2}
\alpha_j\leq d\quad\mbox{for all}\quad j,
\end{equation}
because a rational function of
degree $d$ cannot have a point where the local degree is greater than $d$.

Subject to these two restrictions (\ref{1}) and (\ref{2}), a rational
function with prescribed critical points at $a_j$
of multiplicities $\alpha_j-1$
always exists \cite{G,S}. Thus
{\em for every $a_j$ and every $\alpha_j$ satisfying (\ref{1}) and (\ref{2})
there exist metrics of curvature $1$ on $S$ with angles
$2\pi\alpha_j$ at $a_j$.}

Two rational functions $f_1$ and $f_2$
are called {\em equivalent}
if $f_1=\phi\circ f_2$, where $\phi$ is a fractional-linear
transformation. Equivalent functions have the same critical points
with the same multiplicities. Equivalent functions generate equivalent metrics.

The number of equivalence classes of rational functions with prescribed
critical points and multiplicities is at most
$\K(\alpha_0-1,\ldots,\alpha_{n-1}-1)$,
where $\K$ is the {\em Kostka number} which can be described as follows.
Consider Young diagrams of shape $2\times(d-1)$. They consist of two rows
of length $d-1$ one above another. A semi-standard Young tableau (SSYT)
is a filling of such a diagram with positive integers such that an integer $k$
appears $m_{k-1}$ times, the entries are strictly increasing in the columns
and non-decreasing in the rows. The number of such SSYT's is the
Kostka number $\K(m_0,\ldots,m_{n-1})$.

For a generic choice of the critical points $a_j$, the number
of classes of rational functions is equal to the Kostka number, see
\cite{S,EGSV}.

Suppose now that the points $a_j$ and the corresponding multiplicities
$\alpha_j$
are symmetric with respect to some circle.
We may assume without loss of generality that this circle
is the real line $\R\cup\{\infty\}$.
It may happen that among the rational functions $f$
with these given
critical points and multiplicities none are symmetric. So the resulting
metrics are all asymmetric as well \cite{EG,EG2}.

However, there is a surprising result \cite{EG0,EG,EGSV,MTV}
which was conjectured
by B. and M. Shapiro:
\vspace{.1in}

{\em If all critical points of a rational function lie on a circle,
then the function is equivalent to a function symmetric with respect to this circle. Moreover,
in this case there
are exactly $\K$ equivalence classes of rational functions
with prescribed critical points.}
\vspace{.1in}

It is interesting to find out which of these results
can be extended to the general case of arbitrary positive $\alpha_j$.

\section{The case of two non-integer $\alpha_j$. Condition on the angles.}\label{cond-angles}

In the rest of the paper
we study the case with two non-integer
$\alpha_j$. Some examples of polygons with this property can be found in
\cite{CWX}. We answer the following questions:

a) In the equation (\ref{ode}), for which $\alpha_j$ one can choose
the accessory parameters so that the monodromy group
is conjugate to a subgroup of $PSU(2)$?

b) If $\alpha_j$ satisfy a), how many choices of accessory parameters
satisfying a) are
possible?

c) If in addition all $a_j$ are real, how many choices of real accessory
parameters satisfying a)
are possible?

One cannot have exactly one non-integer $\alpha_j$. Indeed, in this case
the developing map $f$ would have just one branching point on the sphere,
which is impossible by the Monodromy theorem.

Let us consider the case of two non-integer angles.
In this section we obtain
a necessary and sufficient condition on the angles for this case,
that is solve
problem a).

We place the two singularities corresponding to non-integer $\alpha$ at
$a_0=0$ and $a_{n-1}=\infty$,
and let the total angles at these points be $2\pi\alpha_0$ and
$2\pi\alpha_{n-1}$,
where $\alpha_0$ and $\alpha_{n-1}$ are not
integers. Then the developing map has an analytic continuation
in $\C^*$ from which we conclude that the monodromy group must
be a cyclic group
generated by a rotation $z\mapsto ze^{2\pi i\alpha}$,
with some $\alpha\in(0,1)$. This means that $f(z)$ is multiplied by
$e^{2\pi i\alpha}$ when $z$ traverses a simple loop around the origin.
Thus $g(z)=z^{-\alpha}f(z)$ is a single valued function with at most
power growth at $0$ and $\infty$.
Then we have a representation
$f(z)=z^\alpha g(z)$, where $g$ is a rational function.
Then $\alpha_0=|k+\alpha|,\; \alpha_{n-1}=|j+\alpha|,$
where $k$ and $j$ are integers, so
either $\alpha_0-\alpha_{n-1}$ or
$\alpha_0+\alpha_{n-1}$ is an integer.
The angles at the other singular points of the metric
are integer multiples of $2\pi$, and they are the critical points of $f$
other than $0$ and $\infty$.

Let $g=P/Q$ where $P$ and $Q$ are polynomials without common zeros
of degrees $p,q$ respectively.
Let $p_0,q_0$ be the multiplicities
of zeros of $P$ and $Q$ at $0$. Then $\min\{ p_0,q_0\}=0$, because we
may assume that the fraction $P/Q$ is irreducible.

The equation for the critical points of $f$ is the following:
\begin{equation}\label{crit}
z(P'(z)Q(z)-P(z)Q'(z))+\alpha P(z)Q(z)=0.
\end{equation}
Let $2\pi\alpha_0,2\pi\alpha_1\ldots 2\pi\alpha_{n-1}$
be the angles at the conical singularities
at
\newline
$0,a_1,\ldots,a_{n-2},\infty,$
so that $\alpha_j\geq 2$ are integers for $1\leq j\leq n-2.$

Denote
$$\sigma:=\alpha_1+\ldots+\alpha_{n-2}-n+2\geq 0.$$
Then we have the following system of equations
\begin{eqnarray}\label{system}
\alpha_0&=&|p_0-q_0+\alpha|\nonumber\\
\sigma&=&p+q-\max\{ p_0,q_0\}\label{syst}\\
\alpha_\infty&=&|p-q+\alpha|.\nonumber
\end{eqnarray}
The first and the last equations follow immediately from the representation
$f(z)=z^{\alpha}P(z)/Q(z)$ of the developing map.
The second equation holds because the
left hand side of (\ref{crit}) is a polynomial of degree exactly $m+n$,
therefore
the sum of multiplicities of its zeros at $a_j$ for $1\leq j\leq n-2$
must be $p+q-\max\{ p_0,q_0\}$

Solving this system of equations (\ref{syst}) in non-negative integers
satisfying $\min\{ p_0,q_0\}=0$, $p_0\leq p,$ $q_0\leq q$,
we obtain necessary and sufficient conditions
the angles should satisfy, which we state as

\begin{thm}\label{theorem1}Suppose that
$n$ points $a_0,\ldots,a_{n-1}$ on the Riemann
sphere and numbers $\alpha_j>0$, $0\leq j\leq n-1$, are such that
$\alpha_j\geq 2$ are integers for $1\leq j\leq n-2$.

The necessary and sufficient conditions
for the existence of a metric of curvature $1$ on the sphere, with conic
singularities at $a_j$ and angles $2\pi\alpha_j$ are the following:

a) If $\sigma+[\alpha_0]+[\alpha_\infty]$ is even, then $\alpha_0-\alpha_\infty$
is an integer, and
\begin{equation}\label{53}
|[\alpha_0]-[\alpha_\infty]|\leq \sigma.
\end{equation}

b) If $\sigma+[\alpha_0]+[\alpha_\infty]$ is odd, then
$\alpha_0+\alpha_\infty$ is an integer, and
\begin{equation}\label{54}
[\alpha_0]+[\alpha_\infty]+1\leq \sigma.
\end{equation}
\end{thm}

{\em Proof.} Let us first verify
that conditions a)-b)
are necessary and sufficient for the existence of a unique solution
$p,q,p_0,q_0,\alpha$
of the system (\ref{syst}) satisfying
$\min\{ p_0,q_0\}=0,\quad p_0\leq p,\quad q_0\leq q,\quad
\alpha\in(0,1).$

We consider four cases.

Case 1. $p\geq q,\; m_0\geq q_0=0$. Then
\begin{eqnarray*}
p_0&=&[\alpha_0],\\
p+q&=&\sigma+[\alpha_0],\\
p-q&=&[\alpha_\infty].
\end{eqnarray*}
So $\alpha_0-\alpha_\infty$ is an integer, and
$2p=\sigma+[\alpha_0]+[\alpha_\infty]$ is even.
Now
$$0\leq 2(p-p_0)=\sigma-[\alpha_0]+[\alpha_\infty],$$
$$0\leq 2q=(\sigma+[\alpha_0]-[\alpha_\infty]),$$
so we obtain (\ref{53}).
\vspace{.1in}

Case 2. $p\geq q,\; q_0>p_0=0$. Then
\begin{eqnarray*}
q_0&=&[\alpha_0]+1,\\
p+q&=&\sigma+[\alpha_0]+1,\\
p-q&=&[\alpha_\infty].
\end{eqnarray*}
So $\alpha_0+\alpha_\infty$ is an integer, and
$2p=\sigma+[\alpha_0]+[\alpha_\infty]+1$ is even.
We have
$2q=(\sigma+[\alpha_0]-[\alpha_\infty]+1,$
so
$$0\leq 2(q-q_0)=\sigma-[\alpha_0]-[\alpha_\infty]-1.$$
which gives (\ref{54}).
\vspace{.1in}

Case 3. $p<q,\; p_0\geq q_0=0$. Then
\begin{eqnarray*}
p_0&=&[\alpha_0],\\
p+q&=&\sigma+[\alpha_0],\\
q-p&=&[\alpha_\infty]+1.
\end{eqnarray*}
So $\alpha_0+\alpha_\infty$ is an integer, and
$2q=\sigma+1+[\alpha_0]+[\alpha_\infty]$ is even.
We have
$2p=\sigma-1+[\alpha_0]-[\alpha_\infty],$
so
$$0\leq 2(p-p_0)=\sigma-1-[\alpha_0]-[\alpha_\infty],$$
which is (\ref{54}).
\vspace{.1in}

Case 4. $p<1,\; q_0>p_0=0$. Then
\begin{eqnarray*}
q_0&=&[\alpha_0]+1,\\
p+q&=&\sigma+1+[\alpha_0],\\
q-p&=&[\alpha_\infty]+1.
\end{eqnarray*}
So $\alpha_0-\alpha_\infty$ is an integer, and
$2q=\sigma+[\alpha_0]+[\alpha_\infty]+2$ is even.
We have
$0\leq 2p=\sigma+[\alpha_0]-[\alpha_\infty],$
and
$$0\leq 2(q-q_0)=\sigma-[\alpha_0]+[\alpha_\infty],$$
which implies (\ref{53}).

Thus the conditions a) and b) are necessary.

Now we set
$$P(z)=z^{\max\{ p_0,q_0\}}\prod_{j=1}^{n-2}(z-a_j)^{\alpha_j-1}.$$
Second equation in (\ref{syst}) gives $\deg P=p+q$.
Now we consider the equation
\begin{equation}\label{crit2}
z(P^\prime Q-PQ^\prime)+\alpha PQ=R.
\end{equation}
This equation must be solved in polynomials $P,Q$ of degrees $p,q$
having roots of multiplicities $p_0,q_0$ at $0$.
Non-zero polynomials of degree at most $m$ modulo proportionality
can be identified with points of the complex projective space $\P^m$.
The map
\begin{equation}\label{wronsk}
W_\alpha:\P^p\times\P^q\to\P^{p+q},\quad (P,Q)\mapsto z(P'Q-PQ')+\alpha PQ
\end{equation}
is well defined. It is a finite map between compact algebraic varieties,
and it can be represented as a linear projection of the Veronese variety.
Its degree is known: it is equal to
\begin{equation}\label{binom}
\binom{p+q}p.
\end{equation}
Thus the equation (\ref{crit2}) always has a complex solution $(P,Q)$.
The function $f=z^\alpha(P/Q)$ is then a developing map with the
required properties.
So conditions a)-b) are sufficient.
This completes the proof of the first statement.

We will later prove that for generic polynomial $R$ there are exactly
(\ref{binom}) preimages under the map $W_\alpha$.

\section{Case $\alpha=1/2$}

Let us consider a spherical polygon, parametrized by the upper half-plane $H$,
with corners
\begin{equation}\label{positive}
0=a_0<a_1<\ldots<a_{n-2},\quad a_{n-1}=\infty,
\end{equation}
and suppose that the angles $\pi\alpha_j$ at $a_j$
are integer multiples of $\pi$ for
$1\leq j\leq n-2$ while $2\alpha_0$ and $2\alpha_{n-1}$ are {\em odd} integers.
This means that $\alpha=1/2$ in (\ref{wronsk}), and the roots of the polynomial
(\ref{crit}) are non-negative. In intrinsic terms, condition (\ref{positive})
means that the corners with non-integer angles
of our spherical polygon
are adjacent.

Let $Q$ be our spherical polygon, and $Q'$ its mirror image. We paste $Q$
and $Q'$ together identifying the sides between the non-integer corners
isometrically. The result is a polygon $Q^*$ with $2n-2$ corners and all
integer angles. With our half-plane model this procedure can be performed
in the following way. Parametrize our original polygon $Q$ by the first
quadrant, with corners as in (\ref{positive}),
and $Q'$ by the second quadrant with corners at
$$0=a_0>-a_1>-a_2>\ldots>-a_{n-2}>-a_{n-1}=\infty.$$
Then the upper half-plane
will parametrize $Q^*$. Thus the developing map $f^*:H\to\S$ of
the polygon $Q^*$ is a real rational function with all critical points real.
It satisfies
\begin{equation}
\label{sym}
f^*(-\overline{z})=-\overline{f^*(z)},
\end{equation}
so $f^*$ is odd. In the opposite direction, if $f^*$ is an odd real rational
function with with all critical points real, then the restriction of
\begin{equation}\label{corrs}
f(z)=f^*(\sqrt{z})
\end{equation}
to $H$ parametrizes a spherical polygon of
the considered type.

Thus we obtain
\begin{prop}\label{odd-rational}
Equivalence classes of spherical polygons with all angles but two
integer, and adjacent non-integer angles
such that $2\alpha_0$ and $2\alpha_n$ are odd
integers, are in one-to-one correspondence with classes of odd real
rational functions with real critical points.
The correspondence is explicitly given by (\ref{corrs}).

\end{prop}

We recall that equivalence classes of
real rational functions with real critical points are
counted by {\em chord diagrams} \cite{EG0,EG,EGSV}.
Let $g$ be such a function.
Then $g^{-1}(\R\cup{\infty})$ restricted to $H$ is a chord diagram:
it consists of smooth arcs in $H$ with ends on $\R\cup\{\infty\}$,
and these arcs are
disjoint except at the ends. These ends are the
critical points. At a critical point of multiplicity $m$,
$m$ arcs meet, so in our case we have
$$m_0=2\alpha_0-1,\quad m_{n-1}=2\alpha_{n-1}-1,\quad m_j=\alpha_j-1,\quad
1\leq j\leq n-2$$
for the critical points on the ray $0\leq x\leq\infty$.
A rational function is odd iff the diagram is invariant
with respect to $z\mapsto-\overline{z}$, and even number of arcs meet
at $0$ and at $\infty$. We call such chord diagrams {\em odd}.

Without these conditions, the number $\K(n_1,\ldots,n_k)$
of all chord diagrams with vertices of orders $n_1,\ldots,n_k$ is
the Kostka number which was defined in section \ref{all-integers}.
A necessary condition for existence of
such a diagram is that $n_1+\ldots+n_k$ is an even number $2d-2$.
We set $\K(n_1,\ldots,n_k)=0$ if the SSYT with
parameters $n_j$ in the definition of Kostka number
does not exist. We recall that Kostka number does not change
after a permutation of $n_j$.

Now we return to odd chord diagrams.
Let $m_0,m_1,\ldots,m_{n-1}$ be the multiplicities of critical points on
$0\leq x\leq\infty$,
and
\newline
$E(m_0,m_1,\ldots,m_{n-2},m_{n-1})$ the number of chord
diagrams corresponding to odd rational functions
with critical points of multiplicities $m_0,\ldots,m_{n-1}$
on $0\leq x\leq\infty$. We express this quantity in terms of Kostka numbers.
\begin{thm}\label{count}
\begin{equation}\label{kost}
E(m_0,m_1,\ldots,m_{n-2},m_{n-1})=\K(r,m_1,\ldots,m_{n-2},s),
\end{equation}
where positive integers $r$ and $s$
satisfy
\begin{equation}\label{ge}
r+s>m_1+\ldots+m_{n-2},
\end{equation}
and are defined as follows:

If
$\mu:=(m_0+m_{n-1})/2+m_1+\ldots+m_{n-2}$
is even, then $r=m_0/2+k$ and $s=m_{n-1}/2+k$, where $k$ is large
enough to satisfy (\ref{ge}).

If $\mu$ is odd, then $r=|m_0-m_{n-1}|/2+k+1$ and $s=k$,
where $k$ is large enough to satisfy (\ref{ge}).
\end{thm}

{\em Proof.} Let us say that an edge of a chord diagram
is {\em crossing} if its one endpoint
is positive and another negative.
The number $\nu$ of crossing edges satisfies
$\nu\equiv\mu\mod 2.$

If $\nu$ is even, we have
\begin{equation}
\label{nueven}
E(m_0,m_1,\ldots,m_{n-2},m_{n-1})=E(m_0+2,m_1,\ldots,m_{n-2}, m_{n-1}+2).
\end{equation}
Indeed we can establish a bijection between the nets counted
by both sides of (\ref{nueven}) as follows.
If $\nu=0$ we can add two edges connecting $a_0=0$ and $a_{n-1}=\infty$.
If $\nu\neq 0$, we can replace two extreme crossing edges
(closest to $a_0$ and $a_{n-1}$) by the
edges going to $a_0$ and $a_{n-1}$, without changing the other endpoints
of these edges.
Notice that these operations do not change the parity of $\nu$.

If $\nu$ is odd, and $m_{n-1}>0$, we have
\begin{equation}\label{nuodd1}
E(m_0,m_1,\ldots,m_{n-2},m_{n-1})=E(m_0+2,m_1,\ldots,m_{n-2},m_{n-1}-2).
\end{equation}
Indeed, in this case we can shift one crossing edge so that it
becomes an edge ending at $a_0$, the other endpoint unchanged,
and at $a_{n-1}$ to perform the opposite operation: replace one edge
ending at $a_{n-1}$ by a crossing edge.

Similarly if $\nu$ is odd and $m_0>0$, we have
\begin{equation}
\label{nuodd2}
E(m_0,m_1,\ldots,m_{n-2},m_{n-1})=E(m_0-2,m_1,\ldots,m_{n-2}, m_{n-1}+2).
\end{equation}
These operations do not change the number $\nu$.
Finally, if $\nu$ is odd and $m_0=0$, we have
\begin{equation}\label{parnu}
E(0,m_1,\ldots,m_{n-2},m_{n-1})=E(0,m_1,\ldots,m_{n-2}, m_{n-1}+2),
\end{equation}
and similarly, if $\nu$ is odd and $m_{n-1}=0$ then
$$E(m_0,m_1,\ldots,0)=E(m_0+2,m_1,\ldots,m_{n-2},0).$$
and this operation replaces $\nu$ by $\nu-1$, so it switches
the parity of $\nu$.

Finally we notice that if $m_0>2(m_1+\ldots+m_{n-2})$ then there can
be no crossing edges, thus our chord diagram is simply the union
of two symmetrical chord diagrams, and thus
$$E(m_0,\ldots,m_{n-1})=\K(m_0/2,m_1,\ldots,m_{n-2},m_{n-1}/2).$$
The result follows from (\ref{nueven}), (\ref{nuodd1}), (\ref{nuodd2})
and (\ref{parnu}).
\vspace{.1in}

There is no explicit formula for Kostka numbers.
But we need a simple formula in the generic case.
We recall that the number $\K(1,1,\ldots,1)$ of chord diagrams
with $2d$ vertices and one chord ending at each vertex is
$$C_d=\frac{1}{d+1} \binom{\,2d\,}d,$$
the Catalan number.

\begin{prop} For $m$ vertices on the positive ray, one edge ending
at each vertex, we have
$$E_m:=E(0,1,\ldots,1,0)=\binom m{\,[m/2]\,}.$$
\end{prop}

{\em Proof.} In view of Theorem~\ref{count} it is sufficient co count
$\K(a,m_0,\ldots,m_{n-2},b)$ with
\begin{equation}\label{in}
a+b>m_1+\ldots+m_{n-2}.
\end{equation}
Let us define the positive integer $d$ by the formula
$$2(d-1)=a+b+m_1+\ldots+m_{n-2}.$$
Then we are counting Young tableaux of shape $2\times (d-1)$.
In such a tableau, $1$'s must stand in the first row on the left
and $n$'s in the second row on the right. In view of the inequality
(\ref{in}) the part of the tableau with entries $2,\ldots,n-1$
consists of two rows that have no common columns.
the number of ways to fill these rows is the binomial coefficient.

\vspace{.1in}

Now we return to solution of the equation (\ref{crit2}) with generic
real $R$, whose roots are simple and positive, and $\deg R=p+q$.
We conclude that $p=q$ or $q=p+1$, and obtain that the degree of the
map $W_\alpha$ is
$$\binom{p+q}p,$$
according to (\ref{binom}). So all preimages $W_{1/2}^{-1}(R)$
are real. Taking limits of such polynomials we conclude that

\begin{prop}\label{real}
For a real polynomial $R$ with non-negative roots, all solutions
of the equation $W_{1/2}(P,Q)=R$ are real.
\end{prop}

This proposition can be also derived from \cite[Theorem 4.2 (ii)]{MTV2}.

\section{Deformation argument}\label{deformation}

Now we extend the results of the previous paragraph to the case of
arbitrary $\alpha$ and prove the following

\begin{thm}
Let $\alpha_0,\alpha_1,\ldots,\alpha_{n-1},\; n\geq 4,$
be given positive numbers
of which $\alpha_0$ and $\alpha_{n-1}$ are not integers and the rest are
integers, and suppose that conditions of Theorem~\ref{theorem1} are satisfied.

Then for given $a_0,\ldots,a_{n-1}$ there are at most
\begin{equation}
E(2[\alpha_0]-1,\alpha_1-1,\ldots,\alpha_{n-2}-1,2[\alpha_{n-1}]-1)
\end{equation}
equivalence classes of metrics with conic singularities
$a_0,\ldots,a_n$, and angles $\pi\alpha_0,\ldots,\pi\alpha_{n-1}$.
For generic $a_j$ we have equality.

Moreover, we always have equality, when the $a_j$ lie on a circle
in the cyclic order $a_0,\ldots,a_{n-1}$.
\end{thm}

Let us fix some $\alpha\in(0,1)$.
Let $Q$ be a spherical $n$-gon with all angles but two integer multiples of
$\pi$, the two corners with non-integer angles adjacent, and
the non-integer angles satisfy that $2\alpha_0$ and $2\alpha_{n-1}$
are odd integers.

We parametrize $Q$ by the upper half-plane $H$ as before,
with corners as in (\ref{positive}).
Let $f:Q\to \S$ be the developing map. Then all sides of $Q$ are mapped
into two circles, $f(\R_+)\subset C_0$ and $f(\R_-)\subset C_{1/2}$.
We may assume that $C_0$ is the real line and $C_{1/2}$ is the imaginary
line.

Let $\psi:\S\to \S$ be a quasiconformal homeomorphism that
leaves the real line pointwise fixed and maps the imaginary
line onto $C_\alpha=\{ re^{i\alpha}:r\in\R\}$.
Then
\begin{equation}\label{qua0}
g=\psi\circ f,
\end{equation}
is a topologically
holomorphic map $H\to \S$ which is a local homeomorphism in the interior
of $H$ and maps the real line into $\R\cup C_{\alpha}\cup\infty$.
By solving a Beltrami equation, we find a quasiconformal
homeomorphism $\psi:H\to H$ such that
\begin{equation}\label{qua}
g=\psi\circ f\circ\phi
\end{equation}
is
holomorphic. Then $f:H\to \S$ is a developing map of a polygon
$Q_\alpha$.

This correspondence can be inverted: for every spherical polygon
$Q_\alpha$ with all angles but two integers, and adjacent corners with
non-integer angles, and developing map $g$ we can find a polygon $Q$
whose non-integer angle at one vertex has fractional part $1/2$,
such that the developing maps of $Q$ and $Q_\alpha$ are related by
formula (\ref{qua}).

Let
\begin{equation}\label{a}
0<a_1<\ldots<a_{n-2}<\infty
\end{equation}
be the corners of $Q$, and
$0<b_1<\ldots<b_{n-2}<\infty$
be the corners of $Q_\alpha$. We evidently have $a_j=\phi(b_j)$.
It is convenient no normalize $a_1=b_1=1$, then the map $\psi$ which
will
uniquely define $\phi$.
With this normalization, the set of possible configurations $a_j$
forms the interior of a simplex $\Delta_a$ of dimension $n-3$.

Let $N$ be the chord diagram corresponding to $f$ as in the previous section.
It is known \cite{EGSV} that for a fixed chord diagram,
and for each point $x\in\Delta_a$, there exists a unique equivalence
class
of rational functions and thus a map $f$.
To this $f$ we put into correspondence the map $g$ by formula (\ref{qua}),
and obtain a point $y=\Phi(x)\in\Delta_b$.

We claim that the map $\Phi$ is a diffeomorphism
of open simplexes $\Delta_a\to\Delta_b$.

First we prove that it is surjective.
Whenever some subset $a_j,a_{j+1},...,a_{j+k}$ collide, the corresponding subset
$b_j,b_{j+1},\ldots,b_{j+k}$ collides, because $\phi$
and $\phi^{-1}$ are normalized quasiconformal maps
with maximum dilatation that depends only on $\alpha$.
So the faces of $\Delta_a$ correspond
to faces of $\Delta_b$ with the same names.
Thus the degree of $\Phi$ must be equal to one, and the map is surjective.

In particular, this implies that for every polynomial $R$ with all roots
positive and simple the full preimage under the map $W_\alpha$ is real,
the fact which can be alternatively derived from \cite[Theorem 4.3 (ii)]{MTV2}.
It follows that the Jacobi determinant of $W_\alpha$ cannot be zero
over a polynomial $R$ with all non-negative roots (no matter simple or not).

We have that $\Phi:\Delta_a\to\Delta_b$ has degree $1$
and Jacobian in never zero. It follows that the map is a diffeomorphism.
This proves the theorem.
\vspace{.1in}


\vspace{.1in}

{\em A. E. and A. G.: Department of Mathematics, Purdue University,

West Lafayette, IN 47907-2067 USA
\vspace{.1in}

V. T.: Department of Mathematics, IUPUI,

Indianapolis, IN 46202-3216 USA;

St. Petersburg branch of Steklov Mathematical Institute,}

\end{document}